\documentclass[11pt]{article}
\usepackage[ansinew]{inputenc}
\usepackage{array}
\usepackage{color}
\usepackage{amsmath,amsthm,amscd}
\usepackage{amsxtra}
\usepackage{amstext}
\usepackage{amssymb}
\usepackage{latexsym}
\usepackage{amsfonts,cite}
\usepackage{graphics,epstopdf}
\usepackage{epsfig}
\usepackage{amsfonts}
\usepackage{graphicx}
\usepackage{hyperref}
\usepackage{tcolorbox}
\usepackage{array,tabularx,colortbl}
\usepackage[framemethod=TikZ]{mdframed}

\providecommand{\U}[1]{\protect\rule{.1in}{.1in}}

\begin{document}
	\sloppy
	\newtheorem{thm}{Theorem}
	\newtheorem{cor}{Corollary}
	\newtheorem{lem}{Lemma}
	\newtheorem{prop}{Proposition}
	\newtheorem{eg}{Example}
	\newtheorem{defn}{Definition}
	\newtheorem{rem}{Remark}
	\newtheorem{note}{Note}
	\numberwithin{equation}{section}
	\newtheorem{Observation}{Observation}
	\thispagestyle{empty}
	\parindent=0mm
	\begin{center}
		{\Large \textbf{Gauss-Appell polynomials: An umbral calculus approach}}\\ 					
		\vspace{0.20cm}
		
		{\bf Subuhi Khan$^{1}$, Ujair Ahmad$^{2}$, Mehnaz Haneef$^{3}$, Serkan Araci$^{4}$}\\
		\vspace{0.15cm}
		
$^{1,2,3}$Department of Mathematics, Aligarh Muslim University, Aligarh, India \\
$^{4}$Department of Computer Engineering, Hasan Kalyoncu University, Gaziantep, T$\ddot{u}$rkiye\\
\footnote{$^{*}$This work has been done under Junior Research Fellowship (Ref No. 231610072319 dated:29/07/2023) awarded to the second author by University Grand Commission, New Delhi.}
	\footnote{$^{1}$E-mail:~subuhi2006@gmail.com (Subuhi Khan)}
	\footnote{$^{2}$E-mail:~ujairamu1998@gmial.com (Ujair Ahmad)}
\footnote{$^{3}$E-mail:~mehnaz272@gmail.com (Mehnaz Haneef)}
\footnote{$^{4}$E-mail:~serkan.araci@hku.edu.tr (Serkan Araci)(Corresponding author)}
			
	\end{center}
	
	\begin{abstract}
		\noindent
This article aims to reinforce the broad applicability of the umbral approach to address complex mathematical challenges and contribute to various scientific and engineering endeavors. The umbral methods are used to reformulate the theoretical framework of special functions and provide powerful techniques for uncovering new extensions and relationships among these functions. This research article introduces an innovative class of special polynomials, specifically the Gauss-Appell polynomials. The fundamental attributes of this versatile family of special polynomials are outlined, including generating relations, explicit representations, and differential recurrence relations. Certain examples of the particular members that belong to the class of Gauss-Appell polynomials are also considered.

	\end{abstract}
	\parindent=0mm
	\vspace{.25cm}
	
	
	\vspace{0.25cm}
	
	\section{Preliminaries}
Umbral calculus has enhanced profound insights into a wide range of mathematical phenomena. One of its fundamental principles is the notion of umbral composition, which allows for the applications of umbral notions through operations that resemble traditional algebraic composition. The umbral technique allows researchers to unveil hidden relationships, simplify complex calculations, and gain deeper insights into various areas of mathematics. The abstract approach of umbral methods is based on the proficiency of formal symbols, which may initially seem enigmatic, but it has proven itself to be a valuable asset in solving problems and advancing the understanding of fundamental mathematical concepts. Its applications in mathematical physics have evolved considerably in recent years. Notably, the umbral method provides novel approaches for solving various partial differential equations frequently arising in physical problems, see for example \cite{DSML,SLicciardi,our3}.\\

Hypergeometric functions are a remarkable class of mathematical functions that are characterized by their versatility and ability to represent a wide range of mathematical functions, including many other special functions as special cases. They are particularly valuable for solving problems involving complex relationships, such as those encountered in probability theory, combinatorics, number theory, and more.\\

The Gauss hypergeometric function is characterized by the following series \cite{2}:
\begin{equation}\label{haeq6}
	{_2F}_1( a,b;c;x)=\sum_{k=0}^{\infty}\frac{(a)_k(b)_k}{(c)_k}\frac{x^k}{k!}.
\end{equation} 
Very recently in \cite{our3}, the umbral form of the Gauss hypergeometric function is introduced in terms of an exponential function by following relation:
\begin{equation}\label{haeq8}
	{_2F}_1(a,b;c;x)=e^{\hat{\chi} x}\phi_{0},
\end{equation}
where the  umbral operator $\hat{\chi}$ acts on the umbral vacuum $\phi_{0}$ such that
\begin{equation}\label{haeq9}
	\hat{\chi}^k\phi_{0}:=\phi_{k}=\frac{(a)_k(b)_k}{(c)_k}.
\end{equation}
Here
\begin{equation}\label{haeq47}
	e^{k\partial \zeta}(a)_{\zeta}|_{\zeta=0}=(a)_{k}.
\end{equation}
The umbral operator $\hat{\chi}$ is the shift operator
\begin{equation*}
	\hat{\chi}:=e^{\partial\zeta_{1}\;\partial\zeta_{2}\;\partial\zeta_{3}},
\end{equation*}
with $\zeta_{1}$, $\zeta_{2}$, and $\zeta_{3}$ as domain variables of the functions that the operator acts on. The vacuum function is given by the following expression:
\begin{equation*}
	\phi(a,\zeta_{1},b,\zeta_{2},c,\zeta_{3})=\frac{(a)_{\zeta_{1}}(b)_{\zeta_{2}}}{(c)_{\zeta_{3}}}.
\end{equation*}
In conclusion, it follows that
\begin{equation*}
	e^{k\;\partial\zeta_{1}}e^{k\;\partial\zeta_{2}}e^{k\;\partial\zeta_{3}}\;\phi(a,\zeta_{1},b,\zeta_{2},c,\zeta_{3})|_{\zeta_{1}=\zeta_{2}=\zeta_{3}=0}:=\hat{\chi}^k\phi_{0}=\frac{(a)_k(b)_k}{(c)_k}.
\end{equation*}
Appell polynomials \cite{APP} clasp a significant place in the realm of number theory and mathematical analysis because of their remarkable applications in these areas. These polynomials also arise in other diverse scientific disciplines. Characterized by their elegant mathematical structure and remarkable properties, Appell polynomials have proven to be indispensable tools for solving complex differential equations, evaluating integrals, and describing behavior of physical systems.\\	
The generating function of Appell polynomials is given as:
\begin{equation}\label{haeq1}
	A(t)\,e^{xt}=\sum_{n=0}^{\infty}A_{n}(x)\dfrac{t^{n}}{n!},
\end{equation}
where $A(t)$ has expansion:
\begin{equation}\label{haeq2}
	A(t)=\sum_{n=0}^{\infty}A_{n}\dfrac{t^{n}}{n!}.
\end{equation}
According to the classical umbral calculus \cite{BL}, the terms $p_{n}$ of a sequence are formally replaced by the power $\hat{p}^{n}$ of a new variable $\hat{p}$, referred as umbra of the sequence $\{p_{n}\}_{n=0}^{\infty}$. The original sequence is recovered by the umbral notation 
\begin{equation*}
	\hat{p}^{n}\psi_{r}|_{r=0}=\psi_{n+r}|_{r=0}:=p_{n+r}|_{r=0}=p_{n}.
\end{equation*} 
The introduction of an umbra for $\{p_{n}\}$ requires a constitutive equation that reflects the properties of the original sequence. These ideas are illustrated with the umbrae of the Bernoulli, Euler, and Hermite numbers; see, for example \cite{DMV,GE}. \\

We extend this application to give the umbral form of the Appell polynomials as:
	\begin{equation}\label{haeq5}
	A_{n}(x)=(x+\hat{a})^{n}\psi_{0},
\end{equation}
where the umbral operator $\hat{a}$ acts on vacuum function $\psi_{0}$ such that:
\begin{equation}\label{haeq3}
	\hat{a}^{n}\psi_{0}=A_{n},
\end{equation}
		in view of equations \eqref{haeq2} and \eqref{haeq3}, we have
		\begin{equation}\label{haeq4}
			A(t)=e^{\hat{a}t}\psi_{0}.
		\end{equation}
		Consequently, the following umbral form of generating function \eqref{haeq1} is obtained:
		\begin{equation}\label{haeq7}
			e^{(x+\hat{a})t}\psi_{0}=\sum_{n=0}^{\infty}A_{n}(x)\dfrac{t^{n}}{n!}.
		\end{equation}
		Certain members belonging to the class of Appell polynomials are listed in Table $1$.
		\begin{table}[h!]
			\scriptsize
			\caption{\bf Umbral forms for certain members of the Appell family}
			\label{tab:table1}
			\begin{tabular}{|c|p{3.00cm}|p{3.79cm}|p{3.50cm}|p{2.8cm}|}
				\hline
				&&&&\\
				\textbf{S.No.} & \textbf{Special polynomials} &	\centering\textbf{A(t)} & \textbf{Generating function}&\textbf{Umbral form}\\
				\hline
				&&&&\\
				I. & Bernoulli polynomials \cite{R} & $\dfrac{t}{e^{t}-1}=e^{\hat{B}t}\psi_{0}$; $\hat{B}^{n}\psi_{0}=B_{n}$& $e^{(x+\hat{B})t}\psi_{0}	=\sum\limits_{n=0}^{\infty}B_{n}(x)\dfrac{t^{n}}{n!}$&$B_{n}(x)=(x+\hat{B})^{n}\psi_{0}$\\
				\hline
				&&&&\\
				II. & Euler polynomials \cite{R} & $\dfrac{2}{e^{t}-1}=e^{\hat{E}t}\psi_{0}$; $\hat{E}^{n}\psi_{0}=E_{n}$& $e^{(x+\hat{E})t}\psi_{0}=\sum\limits_{n=0}^{\infty}E_{n}(x)\dfrac{t^{n}}{n!}$&$E_{n}(x)=(x+\hat{E})^{n}\psi_{0}$\\
				\hline
				&&&&\\
				III. &Genocchi polynomials \cite{DMS}&$\dfrac{2t}{e^{t}+1}=e^{\hat{G}t}\psi_{0}$; $\hat{G}^{n}\psi_{0}=G_{n}$&
				$e^{(x+\hat{G})t}\psi_{0}=\sum\limits_{n=0}^{\infty}G_{n}(x)\dfrac{t^{n}}{n!}$&$G_{n}(x)=(x+\hat{G})^{n}\psi_{0}$\\
				\hline
			\end{tabular}
		\end{table}\\
		
		The first advancement in introducing hybrid families associated with Appell polynomials was made by Subuhi Khan {\em et al.} \cite{SG}. Since then, certain mixed families linked with Appell polynomials using operational methods have been introduced; see, for example, \cite{SRM, KR}. However, due to certain constraints, the Appell family could not be linked with Gauss hypergeometric functions.\\
		
		This article is first attempt in considering the combination of Gauss hypergeometric functions with Appell family by using umbral procedures.
 In Section $2$, the Gauss-Appell polynomials are introduced followed by derivation of their properties. In Section $3$, the differential equations involving Gauss-Appell polynomials are derived within the umbral formalism. Certain particular examples belonging to the class of Gauss-Appell family are considered in Section $4$. Finally, some concluding remarks are given in Section $5$.

\section{Gauss-Appell polynomials}
In order to develop the Gauss-Appell polynomials (GAP) family, the variable $x$ in the Appell polynomials $A_{n}(x)$ is replaced by $x$ times the umbral operator $\hat{\chi}$ of the Gauss hypergeometric functions. Utilizing umbral approach and denoting the resultant GAP by ${_{_{2}F_{1}}{A}}_{n}(a,b;c;x)$, the following umbral representation is obtained:
\begin{equation}\label{haeq14}
	{_{_{2}F_{1}}{A}}_{n}(a,b;c;x)=A_{n}(x\hat{\chi})\phi_{0},
\end{equation}
where $\phi_{0}$ is the vacuum corresponding to the umbral operator $\hat{\chi}$.\\

Using umbral form \eqref{haeq5} of $A_{n}(x)$ in the right hand side of equation \eqref{haeq14}, we get the following umbral form of the GAP ${_{_{2}F_{1}}{A}}_{n}(a,b;c;x)$:
\begin{equation}\label{haeq15}
		{_{_{2}F_{1}}{A}}_{n}(a,b;c;x)=(x\hat{\chi}+\hat{a})^{n}\phi_{0}\psi_{0}.
\end{equation}
Again replacing the variable $x$ by $x\hat{\chi}$ in generating function \eqref{haeq1} of Appell polynomials $A_{n}(x)$, It follows that
\begin{equation*}
	\sum_{n=0}^{\infty}A_{n}(x\hat{\chi})\dfrac{t^{n}}{n!}\phi_{0}=A(t)\,e^{x\hat{\chi}t}\phi_{0}.
\end{equation*}
In view of umbral form \eqref{haeq14}, the above equation gives the following umbral generating function for the GAP ${_{_{2}F_{1}}{A}}_{n}(a,b;c;x)$:
\begin{equation}\label{haeq18}
	A(t)\,e^{x\hat{\chi}t}\phi_{0}=\sum_{n=0}^{\infty}{_{_{2}F_{1}}{A}}_{n}(a,b;c;x)\dfrac{t^{n}}{n!}.
\end{equation}
Now, using umbral form \eqref{haeq4} of $A(t)$ in equation \eqref{haeq18}, the following equivalent umbral generating function of GAP ${_{_{2}F_{1}}{A}}_{n}(a,b;c;x)$ is obtained:
\begin{equation}\label{haeq17}
e^{(x\hat{\chi}+\hat{a})t}\phi_{0}\psi_{0}=	\sum_{n=0}^{\infty}{_{_{2}F_{1}}{A}}_{n}(a,b;c;x)\dfrac{t^{n}}{n!}.
\end{equation}
Also, from equations \eqref{haeq15} and \eqref{haeq18}, we deduce the following umbral equation
\begin{equation}\label{haeq39}	A(t)\,e^{x\hat{\chi}t}\phi_{0}=\sum_{n=0}^{\infty}(x\hat{\chi}+\hat{a})^{n}\phi_{0}\psi_{0}\dfrac{t^{n}}{n!}.
\end{equation}
Further, umbral definitions \eqref{haeq8} and \eqref{haeq4} are used in equation \eqref{haeq17}, to obtain the following ordinary generating function for the Gauss-Appell polynomials:
\begin{equation}\label{haeq12}
	\sum_{n=0}^{\infty}{_{_{2}F_{1}}{A}}_{n}(a,b;c;x)\dfrac{t^{n}}{n!}=A(t)\,{_{2}F_{1}}(a,b;c;xt).
\end{equation}
	Expanding the right hand side of equation \eqref{haeq15}, it follows that:
	\begin{equation*}
		{_{_{2}F_{1}}{A}}_{n}(a,b;c;x)=\sum_{k=0}^{n}\binom{n}{k}(x\hat{\chi})^{k}(\hat{a})^{n-k}\phi_{0}\psi_{0},
	\end{equation*}
	which in view of umbral actions \eqref{haeq3} and \eqref{haeq9}, produces the following series expansion for the Gauss-Appell polynomials ${_{_{2}F_{1}}{A}}_{n}(a,b;c;x)$:
\begin{equation}\label{haeq19}
	{_{_{2}F_{1}}{A}}_{n}(a,b;c;x)=\sum_{k=0}^{n}\binom{n}{k}\dfrac{x^{k}(a)_{k}(b)_{k}A_{n-k}}{(c)_{k}},
\end{equation}
or equivalently
\begin{equation}\label{haeq20}
	{_{_{2}F_{1}}{A}}_{n}(a,b;c;x)=\sum_{k=0}^{n}\binom{n}{k}\dfrac{x^{n-k}(a)_{n-k}(b)_{n-k}A_{k}}{(c)_{n-k}}.
\end{equation}
	Further, replacing $x$ by $x+y$ in umbral definition \eqref{haeq15} and expanding the right hand side of the resultant equation, it follows that
		\begin{equation*}
			{_{_{2}F_{1}}{A}}_{n}(a,b;c;x+y)=\sum_{k=0}^{n}\binom{n}{k}(x+y)^{k}\hat{\chi}^{k}\hat{a}^{n-k}\phi_{0}\psi_{0},
		\end{equation*}
		which on again applying umbral actions \eqref{haeq3} and \eqref{haeq9} yields the following summation formula:
		\begin{equation}\label{haeq21}
			{_{_{2}F_{1}}{A}}_{n}(a,b;c;x+y)=	\sum_{k=0}^{n}\binom{n}{k}\left(1+\dfrac{x}{y}\right)^{k}\dfrac{(a)_{k}(b)_{k}}{(c)_{k}}y^{k}A_{n-k}.
		\end{equation}
	The action of $mth$ power of the operator $\hat{\chi}$ on the umbral form of Gauss-Appell polynomials ${_{_{2}F_{1}}{A}}_{n}(a,b;c;x)$ is given in the following result:
	\begin{thm}
		For the Gauss-Appell polynomials ${_{_{2}F_{1}}{A}}_{n}(a,b;c;x)$, the following identity holds:
	 \begin{equation}\label{haeq41}
				\hat{\chi}^{m}(x\hat{\chi}+\hat{a})^{n}\phi_{0}\psi_{0}=\dfrac{(a)_{m}(b)_{m}}{(c)_{m}}{_{_{2}F_{1}}{A}}_{n}(a+m,b+m;c+m;x).
			\end{equation}
\begin{proof}
	Applying the operator $\hat{\chi}^{m}$ on the both sides of umbral equation \eqref{haeq39}, we have
	\begin{equation*}
		\sum_{n=0}^{\infty}\hat{\chi}^{m}(x\hat{\chi}+\hat{a})^{n}\phi_{0}\psi_{0}\dfrac{t^{n}}{n!}=\hat{\chi}^{m}A(t)\,e^{x\hat{\chi}t}\phi_{0},
	\end{equation*}
which on expanding the exponential in the right hand side and applying operator \eqref{haeq9} in the resultant gives
	\begin{equation}\label{haeq48}
		\sum_{n=0}^{\infty}\hat{\chi}^{m}(x\hat{\chi}+\hat{a})^{n}\phi_{0}\psi_{0}\dfrac{t^{n}}{n!}=A(t)\sum_{n=0}^{\infty}\dfrac{x^{n}(a)_{n+m}(b)_{n+m}}{(c)_{n+m}}\dfrac{t^{n}}{n!}.
	\end{equation}
Using the identity
				\begin{equation}\label{haeqi}
					(r)_{n+m}=(r)_{m}(r+m)_{n}\;,
			\end{equation}
		in equation \eqref{haeq48}, it takes the form
			\begin{equation*}
				\sum_{n=0}^{\infty}\hat{\chi}^{m}(x\hat{\chi}+\hat{a})^{n}\phi_{0}\psi_{0}\dfrac{t^{n}}{n!}=A(t)\dfrac{(a)_{m}(b)_{m}}{(c)_{m}}\sum_{n=0}^{\infty}\dfrac{x^{n}(a+m)_{n}(b+m)_{n}}{(c+m)_{n}}\dfrac{t^{n}}{n!}.
			\end{equation*}
		Further, making use of equation \eqref{haeq12}, the above equation becomes
	\begin{equation*}
		\sum_{n=0}^{\infty}\hat{\chi}^{m}(x\hat{\chi}+\hat{a})^{n}\phi_{0}\psi_{0}\dfrac{t^{n}}{n!}=\dfrac{(a)_{m}(b)_{m}}{(c)_{m}}\sum_{n=0}^{\infty}{_{_{2}F_{1}}{A}}_{n}(a+m,b+m;c+m;x)\dfrac{t^{n}}{n!},
	\end{equation*}
	which on equating coefficients of like powers of $t$ , yields assertion \eqref{haeq41}.
\end{proof}
\end{thm}
	It is to be noted that for $m=1$, assertion \eqref{haeq41} reduces to the following result:
\begin{equation}\label{haeq16}
	\hat{\chi}(x\hat{\chi}+\hat{a})^{n}\phi_{0}\psi_{0}=\dfrac{ab}{c}{_{_{2}F_{1}}{A}}_{n}(a+1,b+1;c+1;x).
\end{equation} 
Next, we derive a recurrence relation for ${_{_{2}F_{1}}{A}}_{n}(a,b;c;x)$  in the form of following result.
\begin{thm}
	For the Gauss-Appell polynomials ${_{_{2}F_{1}}{A}}_{n}(a,b;c;x)$, the following pure recursion formula holds:
	\begin{equation}\label{haeq22}
		{_{_{2}F_{1}}{A}}_{n+1}(a,b;c;x)=\dfrac{abx}{c}{_{_{2}F_{1}}{A}}_{n}(a+1,b+1;c+1;x)+\sum_{k=0}^{n}\binom{n}{k}\beta_{k}\;{_{_{2}F_{1}}{A}}_{n-k}(a,b;c;x),
	\end{equation}
	where
	\begin{equation}\label{haeq23}
		\dfrac{A^{\prime}(t)}{A(t)}=\sum_{k=0}^{\infty}\beta_{k}\dfrac{t^{k}}{k!}.
	\end{equation}
	\begin{proof}
		Taking derivative with respect to t in equation \eqref{haeq18}, we find
		\begin{equation*}
			\sum_{n=1}^{\infty}{_{_{2}F_{1}}{A}}_{n}(a,b;c;x)\dfrac{nt^{n-1}}{n!}=A(t)x\hat{\chi}e^{x\hat{\chi}t}\phi_{0}+A^{\prime}(t)e^{x\hat{\chi}t}\phi_{0}.
		\end{equation*}
		Using equation \eqref{haeq39} in the first term on the right hand side of above equation and simplifying, it follows that
		\begin{equation*}
			\sum_{n=0}^{\infty}{_{_{2}F_{1}}{A}}_{n+1}(a,b;c;x)\dfrac{t^{n}}{n!}=x\sum_{n=0}^{\infty}\hat{\chi}(x\hat{\chi}+\hat{a})^{n}\phi_{0}\psi_{0}\dfrac{t^{n}}{n!}+\dfrac{A^{\prime}(t)}{A(t)}A(t)e^{x\hat{\chi}t}\phi_{0}.
		\end{equation*}
		On using umbral relation \eqref{haeq16} and equation \eqref{haeq23} in the right hand side of above equation, it follows that
	\begin{align*}
	\sum_{n=0}^{\infty}{_{_{2}F_{1}}{A}}_{n+1}(a,b;c;x)\dfrac{t^{n}}{n!}=\dfrac{abx}{c}&\sum_{n=0}^{\infty}{_{_{2}F_{1}}{A}}_{n}(a+1,b+1;c+1;x)\dfrac{t^{n}}{n!}\\&+\left(\sum_{k=0}^{\infty}\beta_{k}\dfrac{t^{k}}{k!}\right)\left(\sum_{n=0}^{\infty}{_{_{2}F_{1}}{A}}_{n}(a,b;c;x)\dfrac{t^{n}}{n!}\right).
	\end{align*}
Applying the Cauchy product rule, the above equation becomes
			\begin{align}\label{haeq25}
		\notag	\sum_{n=0}^{\infty}{_{_{2}F_{1}}{A}}_{n+1}(a,b;c;x)\dfrac{t^{n}}{n!}=\dfrac{abx}{c}\sum_{n=0}^{\infty}{_{_{2}F_{1}}{A}}_{n}&(a+1,b+1;c+1;x)\dfrac{t^{n}}{n!}\\&+\sum_{n=0}^{\infty}\sum_{k=0}^{n}\binom{n}{k}\beta_{k}\;{_{_{2}F_{1}}{A}}_{n-k}(a,b;c;x)\dfrac{t^{n}}{n!}.
		\end{align}
Equating the coefficients of same power of $t$ on both sides of equation \eqref{haeq25}, we get assertion \eqref{haeq22}.
	\end{proof}
\end{thm}
In earlier results, different aspects of the umbral formalism and its significance in relation to Gauss-Appell polynomials are explored. In the following section, we upraise this foundation by extending the formalism to derive differential equations expressed in terms of the umbral representation of Gauss-Appell polynomials.
	\section{Differential equations}
 Differential equations for hybrid families associated with Appell polynomials have been established, see, for example, \cite{KR, S}. The differential equations and other significant attributes of the hybrid special polynomials illustrate their utility in tackling existing and emerging challenges across various scientific fields. This leads us to use the factorization method \cite{IH} to obtain the differential equations for the Gauss-Appell polynomials. To commence, we revisit some fundamental concepts related to this method.\\
	
	Let $\{g_{n}(x)\}_{n=0}^{\infty}$ be a sequence of polynomials such that $deg(g_{n}(x))=n$, ($n\in\mathbb{N}_{0}:=\{0,1,2,\cdots\}$). The lowering shift operator $\rho^{-}$ and raising shift operator $\rho^{+}$
	satisfy the properties
	\begin{equation}\label{haeq11}
		\rho^{-}_{n}\{g_{n}(x)\}=g_{n-1}(x) \quad and \quad \rho^{+}_{n}\{g_{n}(x)\}=g_{n+1}(x).
	\end{equation}
	Further, a useful property can be deduced by using above operators, that is
	\begin{equation}\label{haeq24}
		(\rho^{-}_{n}\;\rho^{+}_{n})\{g_{n}(x)\}=g_{n}(x).
	\end{equation}
	Prior to delving into the main result, the following differential relations are derived.
	\begin{thm}
		For the Gauss-Appell polynomials ${_{_{2}F_{1}}{A}}_{n}(a,b;c;x)$, the following differential relations involving the umbral operators hold:
		\begin{equation}\label{haeq26}
			D_{x\hat{\chi}}^{m}\;[(x\hat{\chi}+\hat{a})^{n}\phi_{0}\psi_{0}]=\dfrac{n!}{(n-m)!}\;{_{_{2}F_{1}}{A}}_{n-m}(a,b;c;x).
		\end{equation}
		\begin{equation}\label{haeq27}
			D_{\hat{a}}^{m}[(x\hat{\chi}+\hat{a})^{n}\phi_{0}\psi_{0}]=\dfrac{n!}{(n-m)!}\;{_{_{2}F_{1}}{A}}_{n-m}(a,b;c;x).
		\end{equation}
		\begin{proof}
	Consider
			\begin{equation}\label{haeqii}
				\sum_{n=0}^{\infty}(x\hat{\chi}+\hat{a})^{n}\phi_{0}\psi_{0}\dfrac{t^{n}}{n!}=e^{(x\hat{\chi}+\hat{a})t}\phi_{0}\psi_{0}.
			\end{equation}
			Taking derivative with respect to {$ x \hat\chi$} in both sides of above equation, we find
			\begin{equation*}
				\sum_{n=0}^{\infty}D_{x\hat{\chi}}\;(x\hat{\chi}+\hat{a})^{n}\phi_{0}\psi_{0}\dfrac{t^{n}}{n!}=te^{(x\hat{\chi}+\hat{a})t}\phi_{0}\psi_{0},
			\end{equation*}
			which on using umbral definition \eqref{haeq17} becomes
			\begin{equation*}
				\sum_{n=0}^{\infty}D_{x\hat{\chi}}\;(x\hat{\chi}+\hat{a})^{n}\phi_{0}\psi_{0}\dfrac{t^{n}}{n!}=	\sum_{n=0}^{\infty}{_{_{2}F_{1}}{A}}_{n}(a,b;c;x)\dfrac{t^{n+1}}{n!}.
			\end{equation*}
			Equating coefficients of $t^{n}$ on both sides of above relation, it follows that:
			\begin{equation}\label{haeq28}
				D_{x\hat{\chi}}\;(x\hat{\chi}+\hat{a})^{n}\phi_{0}\psi_{0}=n\;{_{_{2}F_{1}}{A}}_{n-1}(a,b;c;x),\quad n\ge 1,
			\end{equation}
			which shows that result \eqref{haeq26} holds true for $m=1$. Let us assume that this result remains true for $m=k$, that is
			\begin{equation}\label{haeq29}
				D_{x\hat{\chi}}^{k}\;[(x\hat{\chi}+\hat{a})^{n}\phi_{0}\psi_{0}]=\dfrac{n!}{(n-k)!}{_{_{2}F_{1}}{A}}_{n-k}(a,b;c;x).
			\end{equation}
			Differentiating equation \eqref{haeq29} with respect to {$x \hat\chi$} and using equation \eqref{haeq28} in the right hand side of the resultant equation, we find
			\begin{equation*}
				D_{x\hat{\chi}}^{k+1}\;[(x\hat{\chi}+\hat{a})^{n}\phi_{0}\psi_{0}]=\dfrac{n!}{\left(n-(k+1)\right)!}{_{_{2}F_{1}}{A}}_{n-(k+1)}(a,b;c;x).
			\end{equation*}
			Hence, result \eqref{haeq26} holds true for $m=k+1$. Consequently, by induction, assertion \eqref{haeq26} is proved for all $m \in \mathbb{N}$.\\
			
			\noindent Next, differentiating equation \eqref{haeqii} w.r.t. $\hat{a}$ and proceeding on the same lines as above, assertion \eqref{haeq27} follows.
		\end{proof}
	\end{thm}
	It is important to observe that
	\begin{equation*}
		D_{x\hat{\chi}}^{m}\;[(x\hat{\chi}+\hat{a})^{n}\phi_{0}\psi_{0}]=D_{\hat{a}}^{m}[(x\hat{\chi}+\hat{a})^{n}\phi_{0}\psi_{0}].
	\end{equation*}
In the next result, the expressions of umbral shift operators for the GAP ${_{_{2}F_{1}}{A}}_{n}(a,b;c;x)$ are obtained.
\begin{lem}
	For the Gauss-Appell polynomials ${_{_{2}F_{1}}{A}}_{n}(a,b;c;x)$, the following umbral shift operators exist:
	\begin{equation}\label{haeq30}
	_{x\hat{\chi}}\rho_{n}^{-}=\dfrac{1}{n}D_{x\hat{\chi}},
	\end{equation}
	\begin{equation}\label{haeq31}
		_{\hat{a}}\rho_{n}^{-}=\dfrac{1}{n}D_{\hat{a}},
	\end{equation}
	\begin{equation}\label{haeq32}
		_{x\hat{\chi}}\rho_{n}^{+}=x\hat{\chi}+\sum_{k=0}^{n}\dfrac{\beta_{k}}{k!}D_{x\hat{\chi}}^{k},
	\end{equation}
	\begin{equation}\label{haeq33}
		_{\hat{a}}\rho_{n}^{+}=x\hat{\chi}+\sum_{k=0}^{n}\dfrac{\beta_{k}}{k!}D_{\hat{a}}^{k}.
	\end{equation}
	\begin{proof}
		Taking $m=1$ in equations \eqref{haeq26} and \eqref{haeq27}, we get assertions \eqref{haeq30} and \eqref{haeq31} respectively.\\
		Making use of relation \eqref{haeq26} in the second term on the right hand side of relation \eqref{haeq22}, we have
		\begin{equation*}
			{_{_{2}F_{1}}{A}}_{n+1}(a,b;c;x)=\dfrac{abx}{c}{_{_{2}F_{1}}{A}}_{n}(a+1,b+1;c+1;x)+\sum_{k=0}^{n}\dfrac{\beta_{k}}{k!}\;D_{x\hat{\chi}}^{k}\;(x\hat{\chi}+\hat{a})^{n}\phi_{0}\psi_{0}.
		\end{equation*}
		In view of equations \eqref{haeq15} and \eqref{haeq16}, the above equation takes the form
		\begin{equation*}
			(x\hat{\chi}+\hat{a})^{n+1}\phi_{0}\psi_{0}=x\hat{\chi}(x\hat{\chi}+\hat{a})^{n}\phi_{0}\psi_{0}+\sum_{k=0}^{n}\dfrac{\beta_{k}}{k!}\;D_{x\hat{\chi}}^{k}\;(x\hat{\chi}+\hat{a})^{n}\phi_{0}\psi_{0},
		\end{equation*}
which on simplification	yields assertion \eqref{haeq32}.\\	
			
		Again, by using relation \eqref{haeq27} in \eqref{haeq22} and proceeding on the same lines as above, assertion \eqref{haeq33} follows.
	\end{proof}
\end{lem}
Finally, the umbral differential and umbral partial differential equations for ${_{_{2}F_{1}}{A}}_{n}(a,b;c;x)$ are derived by proving the following result:
\begin{thm}
The Gauss-Appell polynomials  satisfy the following umbral differential equations:
	\begin{equation}\label{haeq34}
	\left(1+x\hat{\chi}D_{x\hat{\chi}}+\sum_{k=0}^{n}\dfrac{\beta_{k}}{k!}D_{x\hat{\chi}}^{k+1}-n\right)(x\hat{\chi}+\hat{a})^{n}\phi_{0}\psi_{0}=0,
	\end{equation}
	or, equivalently
	\begin{equation}\label{haeq35}
	\left(x\hat{\chi}D_{\hat{a}}+\sum_{k=0}^{n}\dfrac{\beta_{k}}{k!}D_{\hat{a}}^{k+1}-n\right)(x\hat{\chi}+\hat{a})^{n}\phi_{0}\psi_{0}=0.
	\end{equation}
	\begin{proof}
	Using operators \eqref{haeq30} and \eqref{haeq32} in factorization relation \eqref{haeq24}, we find
		\begin{equation*}
				{_{x\hat{\chi}}\rho_{n}^{-}}\; {_{x\hat{\chi}}\rho_{n}^{+}}\{(x\hat{\chi}+\hat{a})^{n}\phi_{0}\psi_{0}\}=(x\hat{\chi}+\hat{a})^{n}\phi_{0}\psi_{0},
		\end{equation*}
		which on using expressions of operators \eqref{haeq30} and \eqref{haeq32} yields assertion \eqref{haeq34}.\\		
		Again, using expressions \eqref{haeq31} and \eqref{haeq33} in the factorization relation
		\begin{equation*}
				{_{\hat{a}}\rho_{n}^{-}}\; {_{\hat{a}}\rho_{n}^{+}}\{(x\hat{\chi}+\hat{a})^{n}\phi_{0}\psi_{0}\}=(x\hat{\chi}+\hat{a})^{n}\phi_{0}\psi_{0},
		\end{equation*}
		assertion \eqref{haeq35} is proved.
	\end{proof}
\end{thm}
\begin{thm}
	For the Gauss-Appell polynomials, the following partial differential equations holds true:
	\begin{equation}\label{haeq36}
		\left(1+x\hat{\chi}D_{x\hat{\chi}}+\sum_{k=0}^{n}\dfrac{\beta_{k}}{k!}D_{x\hat{\chi}}\;D_{\hat{a}}^{k}-n\right)(x\hat{\chi}+\hat{a})^{n}\phi_{0}\psi_{0}=0
	\end{equation}
	and
	\begin{equation}\label{haeq37}
		\left(x\hat{\chi}D_{\hat{a}}+\sum_{k=0}^{n}\dfrac{\beta_{k}}{k!}D_{\hat{a}}\;D_{x\hat{\chi}}^{k}-n\right)(x\hat{\chi}+\hat{a})^{n}\phi_{0}\psi_{0}=0.
	\end{equation}
	\begin{proof}
	Since 
	\begin{equation*}
		D_{x\hat{\chi}}(x\hat{\chi}+\hat{a})^{n}\phi_{0}\psi_{0}=D_{\hat{a}}(x\hat{\chi}+\hat{a})^{n}\phi_{0}\psi_{0}
	\end{equation*}	
	and therefore, in view of operators \eqref{haeq30} and \eqref{haeq33}, factorization relation \eqref{haeq24} takes the form
		\begin{equation*}
			{_{x\hat{\chi}}\rho_{n}^{-}}\;{_{\hat{a}}\rho_{n}^{+}}\{(x\hat{\chi}+\hat{a})^{n}\phi_{0}\psi_{0}\}=(x\hat{\chi}+\hat{a})^{n}\phi_{0}\psi_{0},
		\end{equation*}
		which on using expressions \eqref{haeq30} and \eqref{haeq33} yields assertion \eqref{haeq36}.\\
		
		Again, by using expressions \eqref{haeq31} and \eqref{haeq32} in the factorization relation
		\begin{equation*}
			{_{\hat{a}}\rho_{n}^{-}}\;{_{x\hat{\chi}}\rho_{n}^{+}}\{(x\hat{\chi}+\hat{a})^{n}\phi_{0}\psi_{0}\}=(x\hat{\chi}+\hat{a})^{n}\phi_{0}\psi_{0},
		\end{equation*}
	 assertion \eqref{haeq37} follows.
	\end{proof}
\end{thm}
In the next section, examples of certain members belonging to the class of Gauss-Appell polynomials are considered.

\section{Examples}
By taking suitable expressions of the function $A(t)$, specific members belonging to the Gauss-Appell family are obtained. We consider the following examples:\\

{\bf Example 4.1. ~ Gauss-Bernoulli polynomials}\\
In view of generating relation \eqref{haeq17} and Table $1~(\text{S. No. I})$, the The umbral generating function for the Gauss-Bernoulli polynomials ${_{_{2}F_{1}}{B}}_{n}(a,b;c;x)$ is obtained as:
\begin{equation*}
	e^{(x\hat{\chi}+\hat{B})t}\phi_{0}\psi_{0}=	\sum_{n=0}^{\infty}{_{_{2}F_{1}}{B}}_{n}(a,b;c;x)\dfrac{t^{n}}{n!}
\end{equation*}
and the relevant umbral form is given as:
\begin{equation*}
	{_{_{2}F_{1}}{B}}_{n}(a,b;c;x)=(x\hat{\chi}+\hat{B})^{n}\phi_{0}\psi_{0}.
\end{equation*}
For $A(t)=\dfrac{t}{e^{t}-1}$, the generating relation \eqref{haeq12} gives the following equivalent generating function for the Gauss-Bernoulli polynomials ${_{_{2}F_{1}}{B}}_{n}(a,b;c;x)$:
\begin{equation}\label{haeq43}
	\sum_{n=0}^{\infty}{_{_{2}F_{1}}{B}}_{n}(a,b;c;x)\dfrac{t^{n}}{n!}=\dfrac{t}{e^{t}-1}\,{_{2}F_{1}}(a,b;c;xt).
\end{equation}
The Gauss-Bernoulli polynomials can be explicitly represented as:
\begin{equation}
	{_{_{2}F_{1}}{B}}_{n}(a,b;c;x)\newline=\sum\limits_{k=0}^{n}\binom{n}{k}\dfrac{x^{k}(a)_{k}(b)_{k}B_{n-k}}{(c)_{k}}.
\end{equation}
Here, $B_{n-k}:=B_{n-k}(1)$ are the Bernoulli numbers \cite{R}.\\

 The first five Gauss-Bernoulli polynomials are listed below:\\
 \begin{align*}
 	{_{_{2}F_{1}}{B}}_{0}(a,b;c;x)&=1,\\
 	{_{_{2}F_{1}}{B}}_{1}(a,b;c;x)&=-\frac{1}{2}+\frac{ab}{c}x,\\
 	{_{_{2}F_{1}}{B}}_{2}(a,b;c;x)&=\frac{1}{6}-\frac{ab}{c}x+\frac{(a)_2\,(b)_2}{(c)_2}x^2,\\
 	{_{_{2}F_{1}}{B}}_{3}(a,b;c;x)&=\frac{ab}{2c}x-\frac{3}{2}\frac{(a)_2\,(b)_2}{(c)_2}x^2+\frac{(a)_3\,(b)_3}{(c)_3}x^3\\
 	{_{_{2}F_{1}}{B}}_{4}(a,b;c;x)&=-\frac{1}{30}+\frac{(a)_2\,(b)_2}{(c)_2}x^2-2\frac{(a)_3\,(b)_3}{(c)_3}x^3+\frac{(a)_4\,(b)_4}{(c)_4}x^4.
  \end{align*}
  To visulaize these novel polynomials, we draw graphs for these polynomials for some particular values of $n$ using the software Mathematica. For $a=3, b=1, c=7$, Figures 1 and 2 depict the graphs of ${_{_{2}F_{1}}{B}}_{n}(3,1;7;x)$ for $n=2$ and $3$ respectively.\\
  \begin{center}
  	\includegraphics[width=6cm,scale=2]{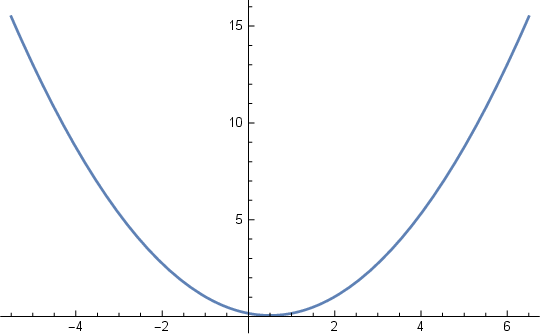}\\{\em \hspace{1cm}Fig. 1 Graph of ${_{_{2}F_{1}}{B}}_{2}(3,1;7;x)$}\\
  	\includegraphics[width=6cm,scale=2]{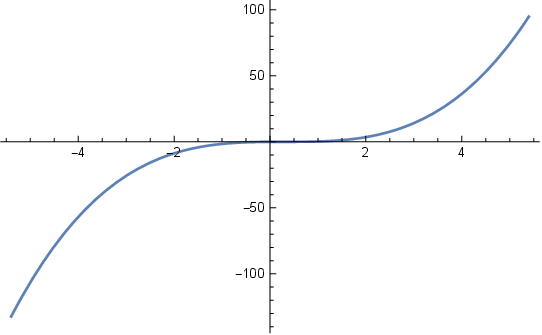}\\
  {\em Fig. 2 Graph of ${_{_{2}F_{1}}{B}}_{3}(3,1;7;x)$}
  \end{center}	
{\bf Example 4.2.~ Gauss-Euler polynomials}\\
The umbral generating function for the  Gauss-Euler polynomials ${_{_{2}F_{1}}{E}}_{n}(a,b;c;x)$ is obtained as:
\begin{equation*}
	e^{(x\hat{\chi}+\hat{E})t}\phi_{0}\psi_{0}=	\sum_{n=0}^{\infty}{_{_{2}F_{1}}{E}}_{n}(a,b;c;x)\dfrac{t^{n}}{n!}
\end{equation*}
and the umbral form is given as:
\begin{equation*}
	{_{_{2}F_{1}}{E}}_{n}(a,b;c;x)=(x\hat{\chi}+\hat{E})^{n}\phi_{0}\psi_{0}.
\end{equation*}
  In view of generating equation \eqref{haeq12}, we obtain the following equivalent generating function for the Gauss-Euler polynomials ${_{_{2}F_{1}}{E}}_{n}(a,b;c;x)$:
 \begin{equation}
 	\sum\limits_{n=0}^{\infty}{_{_{2}F_{1}}{E}}_{n}(a,b;c;x)\dfrac{t^{n}}{n!}=\dfrac{2}{e^{t}+1}{_{2}F_{1}}(a,b;c;xt). 
 \end{equation}
 The explicit series expression for Gauss-Euler polynomials ${_{_{2}F_{1}}{E}}_{n}(a,b;c;x)$ is given as:
 \begin{equation}
 {_{_{2}F_{1}}{E}}_{n}(a,b;c;x)	=\sum\limits_{k=0}^{n}\binom{n}{k}\dfrac{x^{k}(a)_{k}(b)_{k}E_{n-k}}{(c)_{k}},
 \end{equation}
 where $E_{n-k}:=E_{n-k}(1)$ are the Euler numbers \cite{R}.\\
 The expressions of first five Gauss-Euler polynomials are given as:\\
\begin{align*}
	{_{_{2}F_{1}}{E}}_{0}(a,b;c;x)&=1,\\
	{_{_{2}F_{1}}{E}}_{1}(a,b;c;x)&=\frac{ab}{c}x,\\
	{_{_{2}F_{1}}{E}}_{2}(a,b;c;x)&=-1+\frac{(a)_2\,(b)_2}{(c)_2}x^2,\\
	{_{_{2}F_{1}}{E}}_{3}(a,b;c;x)&=-\frac{3ab}{c}x+\frac{(a)_3\,(b)_3}{(c)_3}x^3\\
	{_{_{2}F_{1}}{E}}_{4}(a,b;c;x)&=5-6\frac{(a)_2\,(b)_2}{(c)_2}x^2+\frac{(a)_4\,(b)_4}{(c)_4}x^4.
\end{align*}
We draw the graphs for these polynomials for some particular values of $n$ using the software Mathematica. For $a=3, b=1, c=7$, Figures 3 and 4 depict the graphs of ${_{_{2}F_{1}}{E}}_{n}(3,1;7;x)$ for $n=2$ and $3$ respectively.\\
\begin{center}
	\includegraphics[width=6cm,scale=2]{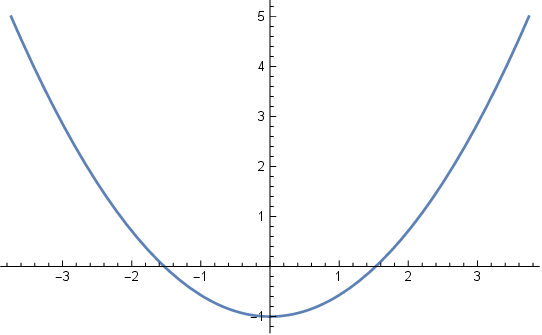}\\{\em \hspace{1cm}Fig. 3 Graph of ${_{_{2}F_{1}}{E}}_{2}(3,1;7;x)$}\\
	\includegraphics[width=6cm,scale=2]{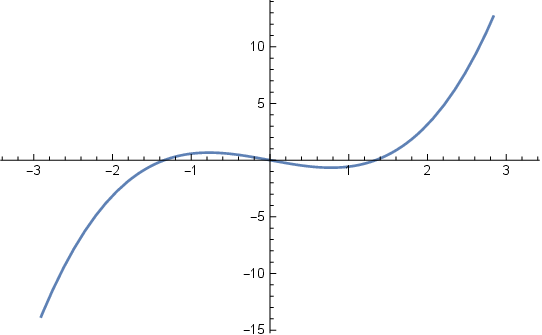}\\
	{\em Fig. 4 Graph of ${_{_{2}F_{1}}{E}}_{3}(3,1;7;x)$}
\end{center}	
{\bf Example 4.3.~Gauss-Genocchi polynomials}\\
 The umbral generating function for the  Gauss-Genocchi polynomials ${_{_{2}F_{1}}{G}}_{n}(a,b;c;x)$ is obtained as:
 \begin{equation*}
 	e^{(x\hat{\chi}+\hat{G})t}\phi_{0}\psi_{0}=	\sum_{n=0}^{\infty}{_{_{2}F_{1}}{G}}_{n}(a,b;c;x)\dfrac{t^{n}}{n!}
 \end{equation*}
 and the umbral form is given as:
 \begin{equation*}
 	{_{_{2}F_{1}}{G}}_{n}(a,b;c;x)=(x\hat{\chi}+\hat{G})^{n}\phi_{0}\psi_{0}.
 \end{equation*}
 In view of the generating equation \eqref{haeq12}, we find the following equivalent generating function for the Gauss-Genocchi polynomials ${_{_{2}F_{1}}{G}}_{n}(a,b;c;x)$:
 \begin{equation}
 		\sum\limits_{n=0}^{\infty}{_{_{2}F_{1}}{G_{n}}}(a,b;c;x)\dfrac{t^{n}}{n!}=\dfrac{2t}{e^{t}+1}{_{2}F_{1}}(a,b;c;xt).
 \end{equation}
 The Gauss-Genocchi polynomials ${_{_{2}F_{1}}{G}}_{n}(a,b;c;x)$ possess the following series expansion formula:
 \begin{equation}
 {_{_{2}F_{1}}{G_{n}}}(a,b;c;x)=\sum\limits_{k=0}^{n}\binom{n}{k}\dfrac{x^{k}(a)_{k}(b)_{k}G_{n-k}}{(c)_{k}},
 \end{equation}
 where $G_{n-k}:=G_{n-k}(1)$ are the Genocchi numbers \cite{RS}.\\
 
 The first five Gauss-Genocchi polynomials are given as:\\
\begin{align*}
	{_{_{2}F_{1}}{G}}_{0}(a,b;c;x)&=0,\\
	{_{_{2}F_{1}}{G}}_{1}(a,b;c;x)&=1,\\
	{_{_{2}F_{1}}{G}}_{2}(a,b;c;x)&=-1+\frac{2ab}{c}x,\\
	{_{_{2}F_{1}}{G}}_{3}(a,b;c;x)&=-\frac{3ab}{c}x+\frac{3(a)_2\,(b)_2}{(c)_2}x^2\\
	{_{_{2}F_{1}}{G}}_{4}(a,b;c;x)&=1-6\frac{(a)_2\,(b)_2}{(c)_2}x^2+4\frac{(a)_3\,(b)_3}{(c)_3}x^3.
\end{align*}
We draw the graphs for these polynomials for some particular values of $n$ using the software Mathematica. For $a=3, b=1, c=7$, Figures 5 and 6 depict the graphs of ${_{_{2}F_{1}}{G}}_{n}(3,1;7;x)$ for $n=2$ and $3$ respectively.\\
\begin{center}
	\includegraphics[width=6cm,scale=2]{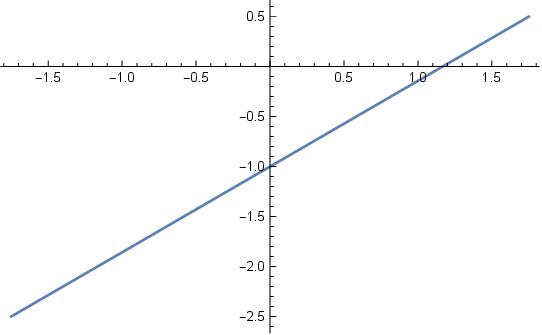}\\{\em \hspace{1cm}Fig. 5 Graph of ${_{_{2}F_{1}}{G}}_{2}(3,1;7;x)$}\\
	\includegraphics[width=6cm,scale=2]{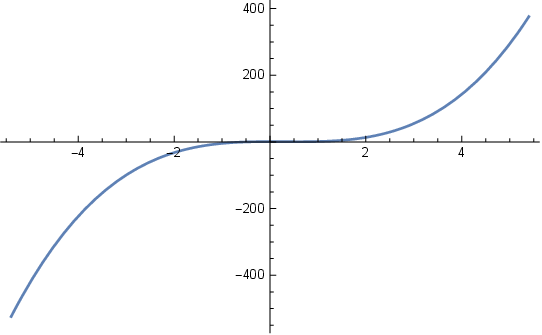}\\
	{\em Fig. 6 Graph of ${_{_{2}F_{1}}{G}}_{3}(3,1;7;x)$}
\end{center}	
Additionally, specific results can also be obtained for members of the Gauss-Appell family to further validate the effectiveness of the umbral formalism.\\

Taking the logarithmic derivative of $A(t)$ w.r.t. $t$ for the Bernoulli polynomials, as listed in Table \thmnumber{1(I)}, yields:
\begin{equation}\label{haeq45}
	\dfrac{A^{\prime}(t)}{A(t)}=-\sum_{k=0}^{\infty}\dfrac{B_{k+1}(1)}{k+1}\dfrac{t^{k}}{k!}.
\end{equation}
Comparison of equations \eqref{haeq45} and \eqref{haeq23} gives
\begin{equation}\label{haeq38}
	\beta_{k}=-\dfrac{B_{k+1}(1)}{k+1}.
\end{equation}

Substituting the expression of $\beta_k$ from equation \eqref{haeq38} in relations \eqref{haeq22} and \eqref{haeq32}--\eqref{haeq37}, we obtain the corresponding results for the Gauss-Bernoulli polynomials ${_{_{2}F_{1}}{B}}_{n}(a,b;c;x)$, which are listed in Table $2$.\\
\begin{table}[h!]
	\centering
	\scriptsize
	\caption{\bf Results for ${_{_{2}F_{1}}{B}}_{n}(a,b;c;x)$}
	\label{tab:table1}
	\begin{tabular}{|p{3.4cm}|p{11.2cm}|}
		\hline
		&\\
		\textbf{Recurrence Relation} & ${_{_{2}F_{1}}{B}}_{n+1}(a,b;c;x)=\dfrac{abx}{c}\;{_{_{2}F_{1}}{B}}_{n}(a+1,b+1;c+1;x)-\sum\limits_{k=0}^{n}\binom{n}{k}\dfrac{B_{k+1}}{k+1}\;{_{_{2}F_{1}}{B}}_{n-k}(a,b;c;x)$ \\
		&\\
		\hline
		&\\
		\textbf{Umbral shift operators} & $	_{x\hat{\chi}}\rho_{n}^{-}=\dfrac{1}{n}D_{x\hat{\chi}}$, \; $_{\hat{B}}\rho_{n}^{-}=\dfrac{1}{n}D_{\hat{B}}$,\newline $	_{x\hat{\chi}}\rho_{n}^{+}=x\hat{\chi}-\sum\limits_{k=0}^{n}\dfrac{B_{k+1}}{(k+1)!}D_{x\hat{\chi}}^{k}$, \;$_{\hat{B}}\rho_{n}^{+}=x\hat{\chi}-\sum\limits_{k=0}^{n}\dfrac{B_{k+1}}{(k+1)!}D_{\hat{B}}^{k}.$\\
		&\\
		\hline
		&\\
		\textbf{Umbral differential equations} & $	\left(1+x\hat{\chi}D_{x\hat{\chi}}-\sum\limits_{k=0}^{n}\dfrac{B_{k+1}}{(k+1)!}D_{x\hat{\chi}}^{k+1}-n\right)(x\hat{\chi}+\hat{B})^{n}\phi_{0}\psi_{0}=0$, \newline
		$\left(x\hat{\chi}D_{\hat{B}}-\sum\limits_{k=0}^{n}\dfrac{B_{k+1}}{(k+1)!}D_{\hat{B}}^{k+1}-n\right)(x\hat{\chi}+\hat{B})^{n}\phi_{0}\psi_{0}=0$.\\
		&\\
		\hline
		&\\
		\textbf{Umbral partial differential equations} & $\left(1+x\hat{\chi}D_{x\hat{\chi}}-\sum\limits_{k=0}^{n}\dfrac{B_{k+1}}{(k+1)!}D_{x\hat{\chi}}\;D_{\hat{B}}^{k}-n\right)(x\hat{\chi}+\hat{B})^{n}\phi_{0}\psi_{0}=0$,\newline $	\left(x\hat{\chi}D_{\hat{B}}-\sum\limits_{k=0}^{n}\dfrac{B_{k+1}}{(k+1)!}D_{\hat{B}}\;D_{x\hat{\chi}}^{k}-n\right)(x\hat{\chi}+\hat{B})^{n}\phi_{0}\psi_{0}=0.$\\
		\hline
	\end{tabular}
\end{table}

\noindent Proceeding on the same lines, analogous results can be derived for the Gauss-Euler and Gauss-Genocchi polynomials.\\
\section{Conclusion}
In the preceding sections, different aspects of umbral formalism are considered. The umbral formalism is particularly effective in simplifying the computational complexities involved. The analysis presented here was focused on the unique properties of this class of functions, especially highlighting their series representations and differential relations. Particular members of the Gauss-Appell family were examined.\\

Further, to extend the legitimacy of the umbral formalism, let us consider the Hermite polynomials $H_{n}(x)$ \cite{2}:
\begin{equation*}
	e^{xt-{t^{2}}/{2}}=\sum_{n=0}^{\infty}\dfrac{H_{n}(x)}{n!}t^{n}.
\end{equation*}
This is a special case of generating function \eqref{haeq1}, as it can be expressed as:
\begin{equation*}
	H(t)e^{xt}=\sum_{n=0}^{\infty}\dfrac{H_{n}(x)}{n!}t^{n},
\end{equation*}
where
\begin{equation*}
	H(t)=e^{-{t^{2}}/{2}}.
\end{equation*}
The umbral form of the Gauss-Hermite polynomials is given as:
\begin{equation*}
{_{_{2}F_{1}}{H}}_{n}(a,b;c;x)=(x\hat{\chi}+\hat{h})^{n}\phi_{0}\zeta_{0}.
\end{equation*}
Here
\begin{equation*}
	\hat{h}^{n}\zeta_{0}=H_{n},
\end{equation*}
where $H_{n}:=H_{n}(0)$ are the Hermite numbers, such that
\begin{equation*}
	H(t)=e^{\hat{h}t}\zeta_{0}.
\end{equation*}
Additionally, certain results can be derived adopting the same approach used throughout this article.\\
 Likewise, the same formalism can be extended to the Laguerre polynomials $L_{n}^{\alpha}(x)$ \cite{2}:
\begin{equation*}
	(1-t)^{\alpha}e^{xt}=\sum_{n=0}^{\infty}(-1)^{n}L_{n}^{\alpha-n}(x)t^{n}.
\end{equation*}
The fact that certain special functions are expressible in terms of Gauss hypergeometric function allows us to visit numerous hybrid polynomials related to the Appell family. Their exploration with a number-theoretic approach will be taken in our subsequent work.\\

The formalism presented in this article may be extended to explore further possibilities including the bivariate special polynomials. In order to provide an illustrations, let us consider the bivariate Hermite-Appell polynomials \cite{SG} defined by the following generating relation:
\begin{equation}\label{haeq49}
	A(t)e^{xt+yt^{2}}=\sum_{n=0}^{\infty}{_{_H}{A}}_{n}(x,y)\dfrac{t^{n}}{n!}.
\end{equation}
Replacing the variable $y$ by $y\hat{\chi}$ in the above equation, we have
\begin{equation}\label{haeq50}
	A(t)e^{xt+y\hat{\chi}t^{2}}\phi_{0}=\sum_{n=0}^{\infty}{_{_H}{A}}_{n}(x,y\hat{\chi})\phi_{0}\dfrac{t^{n}}{n!}.
\end{equation}
Now, following the same procedure adopted in Section $2$ and denoting the resultant bivariate Gauss-Appell polynomials in the right hand side by ${_{_{2}F_{1}}{A}}_{n}^{(a,b;c)}(x,y)$, we get
\begin{equation*}
	\sum_{n=0}^{\infty}{_{_{2}F_{1}}{A}}_{n}^{(a,b;c)}(x,y)\dfrac{t^{n}}{n!}=A(t)e^{xt}{_{_{2}F_{1}}}(a,b;c;yt^{2}).
\end{equation*}
The study of these bivariate polynomials using different aspects including applications in approximation theory may be taken by interested researchers.\\

\noindent {\bf Conflict of Interest:} The authors declare that they have no conflict of interest.


\begin{thebibliography}{99}	 	
	 	\bibitem{2} L. C. Andrews: {\em Special Functions for Engineers and Applied Mathematicians}, Macmillan Publishing Company, New York 1985.
	 	\bibitem{APP} P. Appell: Sur une classe de polyn$\hat{o}$mes, {\em Ann. Sci. $\acute{E}$cole Norm. Sup.}, {\bf 9}(2) 119-144 (1880).
	 	\bibitem{BL} J. Blissard: Theory of generic equations, {\em Quart. J. Pure Appl. Math.} {\bf 4} 279–305 (1861).
		\bibitem{BNR} G. Bretti, P. Natalini and P. E. Ricci: Generalizations of the Bernoulli and Appell polynomials, {\em Abstr. Appl. Anal.} (7) 613-623 (2004).
		\bibitem{DMS} G. Dattoli, M. Migliorati and H. M. Srivastava: Sheffer polynomials, monomiality principle, algebraic methods and the theory of classical polynomials, {\em Math. Comput. Modelling} {\bf 45} 1033-1041 (2007).
	 	\bibitem{DSML} G. Dattoli, Subuhi Khan, Mehnaz Haneef and S. Licciardi: On umbral property of a family of hyperbolic-like functions appearing in magnetic transport problem, {\em Rep. Math. Phys.} {\bf 92}(1) 37-48 (2023). 
	 	\bibitem{our3} G. Dattoli, Subuhi Khan, Mehnaz Haneef and S. Licciardi: Unveiling new perspectives of hypergeometric functions using umbral techniques, {\em Bol. Soc. Mat. Mex.} {\bf 30}(77) 1-21 (2024).	
	 	\bibitem{DMV} A. Dixit, V. Moll and C. Vignat: The Zagier modification of Bernoulli numbers and a
	 	polynomial extension Part I, {\em Ramanujan J.} {\bf 33} 379–422 (2014).
	 	\bibitem{GE} I. Gessel: Applications of the classical umbral calculus, {\em Algebra Universalis} {\bf49} 397–434 (2003).
	 	\bibitem{SRM}Subuhi Khan, N. Raza and M. Ali: Finding mixed families of special polynomials associated with Appell sequences {\em J. Math. Anal. Appl.} {\bf 447}(1) 398–418 (2017).
	 	\bibitem{KR}Subuhi Khan and M. Riyasat: Differential and integral equations for the $2$-iterated Appell polynomials, {\em J. Comput. Appl. Math.} {\bf 306} 116–132 (2016).
	 		\bibitem{SG}Subuhi Khan, G. Yasmin, R. Khan and N. A. M. Hassan: Hermite-based Appell polynomials: Properties and applications, {\em J. Math. Anal. Appl.} {\bf 351} 756–764 (2009).
	 	\bibitem{IH}L. Infeld and T. E. Hull: The factorization method, {\em Rev. Modern Phys.} {\bf 23} 21–68 (1951).
	 	\bibitem{SLicciardi} S. Licciardi and G. Dattoli: {\em Guide to the Umbral Calculus, A Different Mathematical Language}, World Scientific (2022).
	 	\bibitem{R} E. D. Rainville: {\em Special Functions}, The Macmillan Co., New York (1960).
	 \bibitem{RS}J. Riordan and P. R. Stein: Proof of a conjecture on Genocchi numbers, {\em Discrete Math.} {\bf 5} 381–388 (1973).
\bibitem{S}H. M. Srivastava, M. A. $\ddot{O}$zarslan and B. Yilmaz: Some families of differential equations associated with the Hermite-based Appell polynomials and other classes of Hermite-based polynomials, {\em Filomat} {\bf 28}(4) 695–708 (2014).
	 	
 		 \end{thebibliography}
\end{document}